\documentclass[12pt]{article}

\newcommand{\be}{\begin{equation}}
\newcommand{\ee}{\end{equation}}
\newcommand{\bi}[1]{\vspace{-3mm} \bibitem{#1}}

\usepackage{epsfig,amsmath,amssymb,graphics,graphicx} 

\begin{document}

\begin{center}

{\it Communications in Nonlinear Science and Numerical Simulation. \\
Vol.18. No.11. (2013) 2945-2948.} 

\vskip 7mm

{\bf \large No Violation of the Leibniz Rule. No Fractional Derivative.} \\

\vskip 3mm
{\bf \large Vasily E. Tarasov} \\
\vskip 3mm

{\it Skobeltsyn Institute of Nuclear Physics,\\ 
Lomonosov Moscow State University, Moscow 119991, Russia} \\
{E-mail: tarasov@theory.sinp.msu.ru} \\

\begin{abstract}
We demonstrate that a violation of the Leibniz rule is a characteristic property
of derivatives of non-integer orders.
We prove that all fractional derivatives ${\cal D}^{\alpha}$,
which satisfy the Leibniz rule
${\cal D}^{\alpha}(fg)=({\cal D}^{\alpha}f) \, g + f \, ({\cal D}^{\alpha}g)$,
should have the integer order $\alpha=1$, i.e.
fractional derivatives of non-integer orders cannot satisfy the Leibniz rule.
\end{abstract}

\end{center}


\section{Introduction}

Fractional derivatives of non-integer orders \cite{SKM,KST} have 
wide applications in physics and mechanics \cite{CM}-\cite{IJMP2013}.
The tools of fractional derivatives and integrals allow us to investigate 
the behavior of objects and systems that are characterized by power-law non-locality, 
power-law long-term memory or fractal properties.

There are different definitions of fractional derivatives
such as Riemann-Liouville, Riesz, Caputo, Gr\"unwald-Letnikov, Marchaud,
Weyl, Sonin-Letnikov and others \cite{SKM,KST}.
Unfortunately all these fractional derivatives have a lot of unusual properties.
The well-known Leibniz rule ${\cal D}^{\alpha}(fg)=({\cal D}^{\alpha}f)g+f({\cal D}^{\alpha}g)$
is not satisfied for differentiation of non-integer orders \cite{SKM}.
For example, we have the infinite series
\be  \label{GLR}
{\cal D}^{\alpha} (f g) =
\sum^{\infty}_{k=0}  \frac{\Gamma(\alpha+1)}{\Gamma(k+1) \Gamma(\alpha-k+1)} 
( {\cal D}^{\alpha-k} f) \, (D^{k} g)  
\ee
for analytic functions on $[a,b]$ (see  Theorem 15.1 in \cite{SKM}),
where ${\cal D}^{\alpha} $ is the Riemann-Liouville derivative,
$D^{k}$ is derivative of integer order $k$.
Note that the sum is infinite and contains integrals of fractional order for $k>[\alpha]+1$.
Formula (\ref{GLR}) first appeared in the paper by Liouville \cite{Liouv} in 1832.

The unusual properties lead to some difficulties in application of
fractional derivatives in physics and mechanics.
There are some attempts to define new type of fractional derivative
such that the Leibniz rule holds (for example, see \cite{Jumarie1,Jumarie3,Yang}).

In this paper we proof that a violation of the Leibniz rule is one of the main
characteristic properties of fractional derivatives.
We state that linear operator ${\cal D}^{\alpha}$ that can be defined on $C^{2}(U)$, 
where $U \subset \mathbb{R}^1$, such that it satisfied the Leibniz rule
cannot have a non-integer order $\alpha$.  
In other words, a fractional derivative that satisfies the Leibniz rule is not fractional.
It should have integer order.

\section{Hadamard's theorem}

We denote by $C^m(U)$ a space of functions $f(x)$, which are $m$ times continuously differentiable
on $U \subset \mathbb{R}^1$.
Let $D^1_x = d/dx : \ C^m(U) \to C^{m-1}(U)$ be a usual derivative of first order with respect to coordinate $x$.

It is well-known the following Hadamard's theorem \cite{Jet}. 

{\bf Hadamard's Theorem}. 
{\it Any function $f(x) \in C^1(U)$ in a neighborhood $U$ of a point $x_0$ can be represented in the form 
\be \label{HT-1}
f(x)=f(x_0)+  (x-x_0) g(x), 
\ee
where $g(x) \in C^1(U)$.} 

{\bf Proof}. Let us consider the function 
\be
F(t)=f(x_0 +(x - x_0)t). 
\ee
Then $F(0) = f(x_0)$ and $F(1) = f(x)$. 
The Newton--Leibniz formula gives 
\[ F(1)- F(0)= \int^1_0 dt \, (D^1_t F)(t) = \]
\be
= \int^1_0 dt \, (D^1_x f) (x_0+(x-x_0)t) \, (x-x_0) 
=  (x-x_0) \, \int^1_0 dt \, (D^1_x f) (x_0+(x-x_0)t) .
\ee
We define the function
\be
g(x) =\int^1_0 dt \, (D^1_x f) (x_0+(x-x_0)t) .
\ee
As the result, we have proved representation (\ref{HT-1}). \\

\section{Algebraic approach to fractional derivatives}

We consider fractional derivatives ${\cal D}^{\alpha}$ of
non-integer orders $\alpha$ by using an algebraic approach.
Special forms of fractional derivatives are not
important for our consideration.
We take into account the property of linearity and the Leibniz rule only.

For the operator ${\cal D}^{\alpha}$ we will consider the following conditions.

1) $\mathbb{R}$-linearity:
\be \label{cond-1}
{\cal D}^{\alpha}_x  (c_1 f(x) +c_2 g(x)) = c_1 ({\cal D}^{\alpha}_x f(x))  + c_2 \, ({\cal D}^{\alpha}_x  g(x)) ,
\ee
where $c_1$  and $c_2$ are real numbers. 
Note that all known fractional derivatives are linear \cite{SKM,KST}.

2) The Leibniz rule:
\be \label{cond-2}
{\cal D}^{\alpha}_x  (f(x) \, g(x)) = ({\cal D}^{\alpha}_x f(x)) \, g(x) + f(x) \, ({\cal D}^{\alpha}_x  g(x)) .
\ee

3) If the linear operator satisfies the Leibniz rule, then
the action on the unit (and on a constant function) is equal to zero:
\be \label{cond-3}
{\cal D}^{\alpha}_x 1 = 0 .
\ee

Let us proof the following theorem.

{\bf Theorem. ("No violation of the Leibniz rule. No fractional derivative")} \\
{\it If an operator ${\cal D}^{\alpha}_x$ can be applied to functions from $C^{2}(U)$,
where $U \subset \mathbb{R}^1$  be a neighborhood of the point $x_0$,
and conditions (\ref{cond-1}), (\ref{cond-2}) are satisfied,
then the operator ${\cal D}^{\alpha}_x$ is the derivative $D^1_x$ 
of integer (first) order, i.e. it can be represented in the form
\be \label{T-1}
{\cal D}^{\alpha}_x = a(x) \, D^1_x  ,
\ee
where $a(x)$ are functions on $\mathbb{R}^1$.}

{\bf Proof}.

1) Using Hadamard's theorem for the function $g(x)$ 
in the decomposition (\ref{HT-1}), 
the function $f(x)$ for $x \in U$ can be represented in the form 
\be \label{HT-2}
f(x)=f(x_0) +  (x-x_0) g(x_0)+ (x-x_0)^2 \, g_2(x), 
\ee
where $g_2(x) \in C^{2}(U)$,
and $U \subset \mathbb{R}^1$  is a neighborhood of the point $x_0$.
 
Applying to equality (\ref{HT-2}) the operator $D^1_x$ and use $D^1_x f(x_0)=0$, we get  
\be \label{HT-2b}
(D^1_xf)(x)=g(x_0) + 2(x-x_0) \, g_2(x) + (x-x_0)^2 \, (D^1_x g_2)(x) .
\ee
Then
\[ (D^1_x f)(x_0)=g(x_0) . \]
As a result, we have
\be \label{HT-3}
f(x)=f(x_0) +  (x-x_0) (D^1_x f)(x_0)+ (x-x_0)^2 \, g_2(x) . 
\ee

2) Applying to equality (\ref{HT-3}) the operator $ {\cal D}^{\alpha}_x $, we get
\be \label{HT-4}
({\cal D}^{\alpha}_x f)(x)={\cal D}^{\alpha}_x f(x_0) +  {\cal D}^{\alpha}_x \Bigl( (x-x_0) (D^1_x f)(x_0)\Bigr)+ 
{\cal D}^{\alpha}_x \Bigl( (x-x_0)^2 \, g_2(x) \Bigr). 
\ee
The Leibniz rule gives
\[ ({\cal D}^{\alpha}_x f)(x)={\cal D}^{\alpha}_x f(x_0) +  a(x) \, (D^1_x f)(x_0) + 
(x-x_0) \, ({\cal D}^{\alpha}_x \, (D^1_x f)(x_0)+ \]
\be \label{HT-5}
+ 2 a(x) \, (x-x_0) \, g_2(x)   + (x-x_0)^2 \, ({\cal D}^{\alpha}_x  \, g_2) (x) ,
\ee
where we use the notation
\be
a(x)= ( {\cal D}^{\alpha}_x (x-x_0))(x) .
\ee
Then
\be \label{HT-6}
({\cal D}^{\alpha}_x f)(x_0)= { \cal D}^{\alpha}_x f(x_0) +  a(x_0) \, (D^1_x f)(x_0) . 
\ee
As a result, we have
\be
{\cal D}^{\alpha}_x = a(x) \, D^1_x  + b(x) ,
\ee
where we define the function
\be
b(x)= { \cal D}^{\alpha}_x 1 ,
\ee
and we use the $\mathbb{R}$-linearity in the form
\[ { \cal D}^{\alpha}_x f(x_0) = f(x_0) \, ({ \cal D}^{\alpha}_x 1) . \]

3) Using that ${\cal D}^{\alpha}_x 1=0$ for linear operator, which satisfies the Leibniz rule, 
we get $b(x)=0$, ${\cal D}^{\alpha}_x x_0= x_0 {\cal D}^{\alpha}_x 1=0$ and
\be
{\cal D}^{\alpha}_x  = a(x) \, D^1_x , 
\ee
where $a(x) = {\cal D}^{\alpha}_x x$.

As the result, we prove (\ref{T-1}). \\

{\bf Remark.}
In general, the property (\ref{cond-3}) is satisfied not for all type of fractional derivatives.
For example, we have
\[ {\cal D}^{\alpha}_x 1 = \frac{1}{\Gamma (1-\alpha)} x^{- \alpha} \]
for Riemann-Liouville fractional derivative \cite{KST}.
Note that ${\cal D}^{\alpha}_x x$ is not equal to one in general. 
For example,  
\[ {\cal D}^{\alpha}_x x = \frac{1}{\Gamma (2-\alpha)} x^{1- \alpha} \]
for Riemann-Liouville fractional derivative \cite{KST}. \\


Note that this theorem can be proved for multivariable case.

The theorem state that fractional derivative that satisfies the Leibniz rule 
coincides with differentiation of the order equal to one, i.e.
fractional derivatives of non-integer orders cannot satisfy the Leibniz rule.
Unfortunately the Leibniz rule is suggested for some new fractional derivatives 
(the modified Riemann-Liouville derivative 
that is suggested by Jumarie \cite{Jumarie1,Jumarie3}, and
local fractional derivative in the form that is suggested by Yang \cite{Yang}
and some other derivatives).

Linear operators ${\cal D}^{\alpha}_x$ that satisfy the Leibniz rule cannot be considered as 
fractional derivatives of non-integer orders.
Fractional derivatives should be subject to a rule
that is a generalization of the classical Leibniz rule
\[ D^n_x (fg)= \sum^n_{k=0} (D^{n-k} f) \, (D^{k}g) \]
to the case of differentiation and integration of fractional order
(see Section 15. "The generalized Leibniz rule" of \cite{SKM} and references in it). 
It can be assumed with a high degree of reliability that the generalization 
of the Leibniz rule for all types of fractional derivatives 
should be represented by an infinite series in general.
The history of the generalizations of the Leibniz rule
for fractional derivatives,
which is began from the paper \cite{Liouv} in 1868, is described in Section 17 
(Bibliographical Remarks and Additional Information to Chapter 3.) in \cite{KST}.



\end{document}